\documentclass{amsart}

\def\N{\mathbb{N}}

\usepackage{dsfont}
\usepackage{amssymb}

\newtheorem{Theorem}{Theorem}[section]

\newtheorem{Corollary}[Theorem]{Corollary}
\newtheorem{Proposition}[Theorem]{Proposition}
\newtheorem{Rem}[Theorem]{Remark}
\newenvironment{remark}{\begin{Rem}\rm}{\end{Rem}}

\newtheorem{Definition}[Theorem]{Definition}
\newenvironment{definition}{\begin{Definition}\rm}{\end{Definition}}

\newtheorem{Notation}[Theorem]{Notation}
\newenvironment{notation}{\begin{Notation}\rm}{\end{Notation}}

\newtheorem{prf}{\it{Proof}}

\newenvironment{pf}{\begin{prf}\rm}{\hfill$\Box$\end{prf}}
\def\hfl#1{\smash{\mathop{\hbox to 12mm{\lefttarrowfill}}
\limits^{\scriptstyle#1}}}
\def\hfr#1{\smash{\mathop{\hbox to 12mm{\rightarrowfill}}
\limits^{\scriptstyle#1}}}

\begin{document}

\title{Dicritical divisors after S.S.~Abhyankar and I.~Luengo.}
\author{Vincent Cossart and Micka\"el Matusinski $^{ (*)}$}
\thanks{(*) Universit\'e de Versailles Saint-Quentin}

\begin{abstract} In \cite{abh-luengo_dicrit-div}, S.S Abhyankar and I. Luengo introduce a new theory of dicritical divisors in the most general framework. Here we simplify and generalize their results (see Theorems \ref{mainresult} and \ref{mainresult2}).\end{abstract}

\maketitle
\section{Introduction}
The notion of dicritical divisor appeared at the beginning of the 20$^{\textrm{th}}$ century, in the study of isolated singularities of complex planar differential equations  \cite{dulac_pts-dicrit}. Given a germ $\omega$ of a holomorphic differential 1-form singular at the origin $0\in\mathbb{C}^2$, the singularity is called dicritical if there exists an infinite number of irreducible pairwise distinct (germs of) invariant curves passing through $0$. In this case, the resolution of singularities  \cite{seidenberg:reduction} leads  to the following notion (see e.g. \cite{mattei-moussu:holonomie-int-1ere}): a \textbf{dicritical divisor} - if it exists - is an irreducible component of the exceptional divisor which is transverse to the foliation defined by $\omega$. An important example is given by the case where $\omega$ has a first integral which is a meromorphic function  $\frac{f}{g}$, considered in a neighborhood of one of its poles. There the foliation is given by the pencil of curves $\{\lambda f +\mu g=0,\ \lambda,\mu\in \mathbb{C}\}$. This case is related to the Jacobian problem in dimension 2. Indeed, any polynomial map of  $\mathbb{C}^2$ extends to a rational map of  $\mathds{P}_2(\mathbb{C})$ over certain points at infinity (see \cite{le-weber_jacobian-dim2} and Section \ref{sect:jacob} below). 

In connection to the Jacobian problem, S.S. Abhyankar and I. Luengo introduce in \cite{abh-luengo_dicrit-div} an algebraic version of the dicrital divisors, in the most general context. Set a point $p\in \mathds{P}_2(\mathbb{C})$, and consider the corresponding local ring  $R=\mathcal{O}_{\mathds{P}_2(\mathbb{C}),p}$. Pick  $\frac{f}{g}\in QF(R)$ the quotient field of $R$, $\frac{f}{g}$ irreducible. It is well-known that, by a finite sequence of blow-ups of points, one can monomialize the ideal  $(f,g)$ :
$$\mathrm{Spec}\ R \leftarrow X_1\leftarrow\cdots\leftarrow X_\nu.$$
Let  $E=\bigcup_{i}E_i$ be the exceptional divisor in $X_\nu$. For any $i$, we define $\tilde f, \tilde g$ by $f=h \tilde f$, $g=h\tilde g$
where $h=GCD(f,g)$ locally at $x\in E_i$. The couple $(f,g)$ defines a morphism:
$$\begin{array}{llcl}
\phi_{(f,g,i)}:&E_i&\rightarrow& \mathds{P}_1(\mathbb{C})\\
&x&\mapsto& (\tilde{f}(x),\tilde{g}(x)).
\end{array}$$
With this notations, a dicritical divisor is therefore a divisor $E_i$ for which  $\phi_{(f,g,i)}$ is surjective. In other words, through any point of $E_i$ passes the strict transform of a  curve of the pencil. At the origin of their more general definition (\ref{defi:dicrit}), S.S. Abhyankar and I. Luengo make the following key observation: all which precedes is equivalent to supposing that the residu of $\frac{f}{g}$ is transcendental over the residue field of  $\mathcal{O}_{X_\nu,\eta_i}$, the latter being a discrete valuation ring ($\eta_i$ is the generic point $E_i$); the only hypothesis being henceforth that $R$ is a 2 dimensional regular local ring.

\section{Definitions and preleminary results.}
\begin{notation}
From now on, we will use the following notations. Let $R$ be a noetherian regular local ring of dimension 2. We denote by $QF(R)$ its quotient field, $\mathfrak{m}$ its maximal ideal and $K:=R/\mathfrak{m}$ its residue field. We consider also a discrete valuation ring $R_v$ which dominates $R$ and such that $QF(R_v)=QF(R)$. In the other words, $R_v$ is a \textbf{prime divisor} of $R$ in the sense of \cite{abh:val-cent-local-domain}. We denote by $v:QF(R)\rightarrow\mathbb{Z}\cup\{\infty\}$ the corresponding valuation and $\mathfrak{m}_v$ the maximal ideal. The residue map is denoted by  $\mathrm{Res}_v:R_v\rightarrow K_v$, where $K_v:=R_v/\mathfrak{m}_v$ and $K$ is identified with $R/(R\cap \mathfrak {m}_v)$. Note that $\textrm{trdeg}_KK_v=1$. Given a regular system of parameters $(x,y)$ of $\mathfrak{m}$, such a valuation $v$ is said to be (algebraically) \textbf{monomial} with respect to $(x,y)$ if for any polynomial expansion $P(x,y)=\sum_{a,b} \lambda_{a,b}x^ay^b$ with $\lambda_{a,b}\in R\setminus\mathfrak{m}$, one has $v(P)=\min\{av(x)+bv(y)\ |\ \lambda_{a,b}\neq 0\}$ \cite[Definition 3.22]{teissier:valuations}.
\end{notation} 

\begin{definition}\label{defi:dicrit}
Let  $z\in QF(R)$, $z\neq 0$. We call \textbf{dicritical divisor} of $z$ any prime divisor $R_v$ of $R$ such that $z\in R_v$ and $\mathrm{Res}_v(z)$ is transcendantal over $K$.
\end{definition}

We will use the following results and notations \cite[Definition 3, Proposition 3]{abh:val-cent-local-domain} adapted to our context.
\begin{definition}
Let $X\rightarrow \textrm{Spec} (R)$ be the blow-up of $ \textrm{Spec} (R)$ along $\mathfrak{m}$  \cite[Definition p. 163]{hartshorne_alg-geom}. Let $x$ be the center of $v$ over $X$ \cite[Theorem 4.7 p. 101]{hartshorne_alg-geom}. $\hat{R}=\mathcal{O}_{X,x}$ is called  \textbf{blow-up} (or quadratic transform) of $R$ along $v$.\\
More simply, let $(x,y)$ be a regular system of parameters of $\mathfrak{m}$. Suppose for instance that  $v(x)\leq v(y)$. We denote $S:=R[\frac{y}{x}]\cap \mathfrak{m}_v$. We have $\hat{R}:=R[\frac{y}{x}]_{S}$. Then the valuation $v$ has center $S.\hat{R}$.\\
By induction, we call \textbf{sequence of blow-ups} of $R$ along $v$ the sequence $R=R_0\subsetneq \cdots\subsetneq R_i\subsetneq R_{i+1}\subsetneq\cdots$ where for any $i\in\mathbb{N}$, $R_{i+1}$ is the blow-up of  $R_i$ along $v$.
\end{definition}
\begin{Proposition}[Abhyankar]\label{propo:desing}
Let  $(R_j)_{j\in\mathbb{N}}$ be the sequence of blow-ups of $R$ along $v$. There is a unique $\nu\in\mathbb{N}$ such that for any $j,j'\in\mathbb{N}$ with $j\leq \nu<j'$, we have $R_j\neq R_v=R_{j'}$ with $\dim(R_j)=2>1=\dim(R_{j'})$. Moreover, $R_v$ is a prime divisor of $R_\nu$ with $K_v$ pure  transcendental  extension of $K_\nu :=R_\nu/\mathfrak{m}_\nu$ of degree 1.
\end{Proposition}
\begin{remark}\label{rem:eclts}
Let $(x_\nu,y_\nu)$ be a regular system of parameters of  $\mathfrak{m}_\nu$. Then the valuation $v$ is the $\mathfrak{m}_\nu$-adic valuation, which is of course  monomial with respect to $(x_\nu,y_\nu)$, and $R_v=R_\nu[\frac{y_\nu}{x_\nu}]_{(x_\nu)}$. Besides, if we denote  $K_v=K'(t)$ where $K'=R_{\nu}/\mathfrak{m}_{\nu}$ is the relative algebraic closure of $K$ in $K_v$ and $t\in K_v$ is transcendental over $K$, then $[K':K]<\infty$.
\end{remark}

\section{The main theorems.}
The following result is the main theorem of   \cite{abh-luengo_dicrit-div}:

\begin{Theorem}\label{mainresult}
Let $z\in QF(R)$, $z\neq 0$. Let $R_v$ be a dicritical divisor of $z$. Suppose that there is  $x\in \mathfrak{m}\setminus\mathfrak{m}^2$ and $m\in\mathbb{N}$ such that $zx^m\in R$. Then the element $t$ of (\ref{rem:eclts}) can be chosen so that  $\mathrm{Res}_v(z)\in K'[t]$.
\end{Theorem}
\begin{pf}
We proceed by induction on $\nu$ which is finite by (\ref{propo:desing}).

\noindent \underline{Case $\nu=0$}. In this case, the valuation $v$ is the  $\mathfrak{ m}$-adic valuation. Let $(x,y)$ be a regular system of parameters of  $\mathfrak{ m}$. By hypothesis, $z=\displaystyle\frac{f}{x^m}$ and $v(f)=m$. Therefore we write $f$ as:
$$f=\sum_{a+b=m,a,b \in \N}
\lambda_{a,b}x^ay^b,\ \lambda_{a,b}\in R\textrm{, with }\lambda_{a,b}\in R\setminus\mathfrak{m}\textrm{ for some }a,b. \eqno(1)$$
After blowing-up, in $R[\frac{y}{ x}]$,  $v$ is the $x$-adic valuation with valuation ring $R[\frac{y}{ x}]_{(x)}$. So we have:
$$z=\sum_{a+b= m,a,b \in \N}
\lambda_{a,b}x^{a+b-m}\left(\frac{y}{ x}\right)^b=\sum_{a+b= m,a,b \in \N}
\lambda_{a,b}\left(\frac{y}{ x}\right)^b,\ \lambda_{a,b}\in R.$$
Since $\mathrm{Res}(z)$ is transcendental over $K$, there is at least one  $\lambda_{a,b}\in R\setminus\mathfrak{m}$ with $\mathrm{Res}(\lambda_{a,b})\neq 0$ and $b>0$. So we obtain:
$$\mathrm{Res}(z)=\sum_{a+b= m,a,b \in \N}
\mathrm{Res}(\lambda_{a,b})t^b,\ \ \ t:=\mathrm{Res}\left(\frac{y}{x}\right).$$

\noindent \underline{Case $\nu\geq 1$}. In $R=R_0$, we consider the following dichotomy: either $v(y)\geq v(x)$ or $v(y)<v(x)$. 

\noindent Suppose that $v(y)\geq v(x)$. Then $R_1$ is the localisation of   $R\left[\frac{y}{x}\right]$ at the center of $v$. In $R_1$, we have $z=\frac{f_1}{x^{m-o(f)}}$ where $o(f)$ is the $\mathfrak{ m}$-adic order of $f$ and $f_1\in R_1\subset R_v$ is the strict transform of $f$. Since $v(z)=0$ and $v(f_1)\geq 0$, $v(x)\geq 0$, we have $m\geq o(f)$. The hypotheses of the theorem hold in $R_1$: by induction on $\nu$, we obtain the desired result.

\noindent Suppose now that $v(y)<v(x)$. There are two subcases. Either there exists $i,\ 1\leq i \leq \nu$, such that, in the sequence
$$R=R_0 \subset R_1 \subset \cdots R_i \subset \cdots \subset R_{\nu},$$
the inverse image of  $x^m$ in $R_i$ has only one component. Then the hypotheses of the theorem hold for $R_i$: by induction on $\nu$, we obtain the desired result.

Or there is no such $i$. Then the center of $v$ is always at the origin of  one of the two usual affine charts of the blow-ups. The valuation $v$ is monomial defined by 
$$v(x)=\alpha,\ v(y)=\beta,\ \ \ \ \alpha,\beta \in\N\ \mathrm{with}\ \alpha>\beta.  \eqno(2)$$
Since  $\nu\geq 1$ and since we are at the origin of a chart in $R_1$, we have  $R_1=R\left[\frac{x}{y}\right]_{(\frac{x}{y},y)}$. With the notations of  (1), we have
$$\begin{array}{c}
f=\sum_{a\alpha+b\beta \geq m\alpha,a,b \in \N}
\lambda_{a,b}x^ay^b,\ \  \lambda_{a,b}\in R,\\ 
\mathrm{with\ } \lambda_{a,b}\in R\setminus\mathfrak{m}\mathrm{\ for\ at\ least\ one\ couple\ } (a,b)\mathrm{\ such\ that\ }a\alpha+b\beta=m\alpha.
\end{array}$$
By (2), since $\alpha>\beta$, if we have $a\alpha+b\beta - m\alpha=(a+b-m)\alpha+b(\beta-\alpha)=0$, then $(a+b-m)\alpha\geq 0$. So we always have $a+b\geq m$ in the preceding sum.
%De plus, puisque $\mathrm{Res}(z)$ est transcendant sur $K$, il existe au moins un $\lambda_{a,b}\in R\setminus\mathfrak{m}$ avec $\mathrm{Res}(\lambda_{a,b})\not=0 $ et $b>0$. 
So we have:
$$z=\frac{\sum_{a\alpha+b\beta\geq m\alpha}
\lambda_{a,b} (\frac{x}{y})^{a}y^{a+b-m}}{ (\frac{x}{y})^{m}},\ \lambda_{a,b}\in R,$$
$$\mathrm{\ with\ }\lambda_{a,b}\in R\setminus\mathfrak{m}\mathrm{\ for\ at\ least\ one\ couple\ } (a,b)\mathrm{\ such\ that\ }a\alpha+b\beta=m\alpha. $$
The hypotheses of the theorem hold in $R_1$ for such a $z$: by induction on $\nu$, we obtain the desired result.
\end{pf}

The following result generalises Theorem \ref{mainresult} and \cite[Remark (7.4) (II)]{abh-luengo_dicrit-div}.
%Afin d'énoncer le résultat suivant qui généralise le théorème précédent et la remarque (7.4) (II) de ??Abhyankar-Luengo, nous rappelons la définition suivante.
%\begin{Definition}
%Etant donnés $f\in R$, $f\neq 0$, et un système régulier de paramètres $(x,y)$ de $\mathfrak{m}$, on écrit le développement de Taylor correspondant de $f$:
%$$f=\sum_{a,b\in \N}
%c_{a,b}x^ay^b,\ c_{a,b}\in coefficient\ set \eqno(4)$$
% Le \textbf{polygon de Newton} $\Delta(f)$ de $f$ est de l'enveloppe convexe des quadrants  $(a,b)+\mathbb{R}^2$ où $(a,b)\in \mathrm{support}(f)$.
%\end{Definition}
\begin{Theorem}\label{mainresult2}
 Let $z\in QF(R)$, $z\neq 0$. Let $R_v$ be a dicritical divisor of $z$. Suppose that there exist a regular system of parameters $(x,y)$ of $\mathfrak{m}$ and   $a_0,b_0\in\mathbb{N}$ such that  $f=zx^{a_0}y^{b_0}\in R$.
\begin{enumerate}
    \item If $v$ is not monomial with respect to $(x,y)$, then the element $t$ of (\ref{rem:eclts}) can be chosen so that $\mathrm{Res}_v(z)\in K'[t]$.
\item If $v$ is monomial with respect to $(x,y)$, we denote $v(x)=\alpha$, $v(y)=\beta$ ($\alpha,\beta\in\mathbb{N}^*$), $\gamma:=a_0\alpha+b_0\beta=v(f)$ and:%On considère le point $s_0=(k_0,l_0)$, l'entier naturel $\gamma:=k_0\alpha+l_0\beta$, et la droite $D$ d'équation $\alpha k+\beta l=\gamma$ dans le plan $\mathbb{R}^2$ de coordonnées $(k,l)$.  telle que $D\cap\Delta(f)$ soit un côté noté $[s_1,s_2]$ ou un sommet $\{s_1\}$ de $\Delta(f)$. 
$$\begin{array}{c}
 f=\sum_{a\alpha+b\beta \geq \gamma}
\lambda_{a,b}x^ay^b,\ \ \textrm{ with }\lambda_{a,b}\in R \\
\textrm{and }\ B_f:=\{b\in\mathbb{N}\ |\ \lambda_{a,b}\in R\setminus\mathfrak{m},\ a\alpha+b\beta = \gamma\}\neq \emptyset.
\end{array}$$ 
\begin{enumerate}
    \item If $\mathrm{Card}(B_f)\geq 2$, then the element $t$ of (\ref{rem:eclts}) can be chosen so that $\mathrm{Res}_v(z)\in K'[t]$ if and only if  $b_0 \leq \min (B_f)$ or $b_0\geq \max(B_f)$.
\item If $B_f=\{b_1\}$, then we have $b_0\neq b_1$ and the element $t$ of (\ref{rem:eclts}) can be chosen so that $\mathrm{Res}_v(z)\in K'[t]$.
\end{enumerate}  
%le point de coordonnées $(k_0,l_0)$ n'appartient pas $]s_1,s_2[$ ;
%    \item si $D\cap\Delta(f)=\{s_0\}$, alors on a $s_0\neq s_1$ et l'élément $t$ peut être choisi de sorte que $\mathrm{Res}_v(z)\in K'[t]$.
\end{enumerate}
\end{Theorem}
\begin{remark}\label{rem:polyg-newton}
In the case (2), one can translate the condition on $b_0$ in terms of the Newton polygon associated to $f$ relatively to  $(x,y)$. Denote by $D$ the line of equation $\alpha a+\beta b=\gamma$ in the plan $(a,b)\in\mathbb{R}^2$ and by $s_1$, respectively $s_2$, its point with coordinates $(a_1,b_1)$ where $b_1=\min( B_f)$, respectively $(a_2,b_2)$ where $b_2=\max(B_f)$. So $s_0:=(a_0,b_0)\in D$ and the segment $[s_1,s_2]$ is an edge (possibly reduced to a vertex) of the Newton polygon of $f$. The condition   $b_0 \leq \min (B_f)$ or $b_0\geq \max(B_f)$ is equivalent to supposing that  $s_0\notin ]s_1,s_2[$.
\end{remark}
\begin{remark}\footnote{We thank the referee of the J. of Algebra who pointed out the problem of the unicity of $B_f$.}
Note that $$\tilde{f}:=\sum_{(a,b), b\in B_f} \textrm{Res}(\lambda_{a,b})U^aV^b \in K[U,V]\eqno{(E)}$$
may be seen as $\textrm{In}_v(f)\in \textrm{gr}_v R=K[U,V]$ where  $\textrm{In}_v(x)=U$ and  $\textrm{In}_v(y)=V$. By definition in \cite[p.108]{spivakovsky_funct-fields-surf}, $\textrm{gr}_v R:=\oplus_{\rho\in\mathbb{N}}I_\rho/I_{\rho_+}$ where $I_\rho:=\{w\in R\ |\ v(w)\geq \rho\}$ and $I_{\rho_+}:=\{w\in R\ |\ v(w)> \rho\}$. The relation ($E$) is proven in \cite[Remark 3.23 (2)]{teissier:valuations}. See also \cite[Definition 2.9 and Remark 2.10]{hironaka:char} with $\Delta=\{(a,b)\ |\ a\alpha+b\beta\geq 1\}$, which inspired M. Spivakovsky, B. Teissier and many others.

If $a_0.b_0\neq 0$, since $z$ and $f$ are fixed, then $x,y$ are fixed up to multiplication by invertibles. So $U,V$ are fixed up to multiplication by a scalar. Therefore $B_f$ is fixed.

If $a_0=0$ or $b_0=0$, by symmetry between $x$ and $y$ (see Remark \ref{rem:polyg-newton}), we may assume that $b_0=0$. In this case, our condition (2)(a) is trivially verified, even if $B_f$ may not be uniquely defined anymore (if $\alpha | \beta$, we can replace $y$ by $y+\lambda_{a,b}x^{\beta/\alpha}$, i.e. $V$ by $V+\textrm{Res}(\lambda_{a,b})U^{\beta/\alpha}$ which modifies $B_f$).
\end{remark}
\begin{pf}
(1) Since the valuation is not monomial with respect to $x,y$, there is an index $i\in\{1,\ldots,\nu\}$ such that, in $R_i$, we write $x^ky^l=x_i^mu$ with $x_i$ parameter of $\mathfrak{m}_i$ and $u\in R_i$ invertible. Then we are reduced to the hypotheses of Theorem \ref{mainresult}.

\noindent(2) Suppose now that the valuation is monomial with respect to $x,y$. We consider the ring $R_\nu$, with parameters  $x_\nu,y_\nu$, for which the valuation $v$ is $\mathfrak{m}$-adic. We denote $x=x_\nu^{k_1}y_\nu^{k_2}$ and $y=x_\nu^{l_1}y_\nu^{l_2}$ with $k_1 l_2 -k_2l_1=1.$ The exponents $(c,d)$ of the monomials  $x_\nu^cy_\nu^d$ in $R_\nu$ are obtained from the exponants $(a,b)$ of the corresponding monomials $x^ay^b$ by application of a special linear matrix (with the notations of Remark \ref{rem:polyg-newton}, it is the planar linear tranformation changing the line $D$ into $\tilde{D}: c+d=\gamma$):
$$ \left( \begin{array}{c}c\\d\end{array}\right)=A.\left( \begin{array}
{c}a\\b\end{array}\right) \textrm{ where } A=\left( \begin{array}{cc}
k_1 & l_1\\k_2 & l_2\end{array}\right)\textrm{ with } \det(A)=k_1 l_2 -k_2l_1= 1.$$
So we obtain $\mathrm{ord}_{x_\nu,y_\nu}(f)=c_0+d_0=\gamma$ and:
$$ z=\displaystyle\frac{\sum_{c+d\geq c_0+d_0}\lambda_{a,b} x_\nu^cy_\nu^d}{x_\nu^{c_0}y_\nu^{d_0}}\ \textrm{ ; }\ \left( \begin{array}{c}c\\d\end{array}\right)=A.\left( \begin{array}
{c}a\\b\end{array}\right).$$
(2)(a) A linear map preserves barycenters, so  $b_0 \leq \min (B_f)$, respectively $b_0\geq \max(B_f)$, if and only if $d_0\leq \min(\tilde{B}_f)$, respectively $d_0\geq \max(\tilde{B}_f)$, where $\tilde{B}_f=\{d=k_2a+l_2b\ |\ b\in B_f\}$. In the first case, we denote the change of coordinates of the last blow-up by $(x_\nu,y_\nu)\mapsto (x_\nu, \frac{x_\nu}{y_\nu})$ with $\frac{x_\nu}{y_\nu}\in R\setminus\mathfrak{m}$. Setting $t:=\mathrm{Res}\left(\frac{x_\nu}{y_\nu}\right)$, we compute as desired:
$$\begin{array}{lcll}
\mathrm{Res}(z)&=&\displaystyle\frac{\sum_{c+d=\gamma}
\mathrm{Res}(\lambda_{a,b})t^d}{t^{d_0}},&\textrm{with } (c,d)=(k_1 a+l_1b,k_2 a+l_2b);\\
&=&\sum_{c+d= \gamma}
\mathrm{Res}(\lambda_{a,b})t^{d-d_0}\in K[t]\setminus K,& d\geq d_0.
\end{array}$$ 
In the second case, we make the other change of coordinates and we obtain also $\mathrm{Res}(z)\in K[t]\setminus K$ with  $t:=\mathrm{Res}(\frac{y_\nu}{x_\nu})$. On the other hand, if $\min (B_f)<b_0 <\max(B_f)$ (which implies that $\max(B_f)-\min(B_f)\geq 2$), with for instance $t:=\mathrm{Res}(\frac{y_\nu}{x_\nu})$, we obtain:
$$\mathrm{Res}(z)=
\mathrm{Res}(\lambda_{a_1,b_1})t^{d_1-d_0}+\cdots+
\mathrm{Res}(\lambda_{a_2,b_2})t^{d_2-d_0}, \ \textrm{ with } d_1-d_0<0<d_2-d_0.$$ 
We note like \cite[p.1]{abh-luengo_dicrit-div} that $t'$ is another generator of  $K_v=K(t)$ over $K$ if and only if there exists $\rho_1,\rho_2,\theta_1,\theta_2\in K$ such that $t:=\frac{\rho_1t'+\rho_2}{\theta_1t'+\theta_2}$ and  $\rho_1\theta_2-\rho_2\theta_1\neq 0$. Thus we can write $\mathrm{Res}(z)$ as follows: 
$$\mathrm{Res}(z)=\frac{\mathrm{Res}(\lambda_{a_1,b_1}) (\theta_1t'+\theta_2)^{d_2-d_1}+\cdots+
\mathrm{Res}(\lambda_{a_2,b_2})(\rho_1t'+\rho_2)^{d_2-d_1}}{(\theta_1t'+ \theta_2)^{d_0-d_1} (\rho_1t'+\rho_2)^{d_2-d_0}}.$$ 
In the case where $\theta_1\rho_1=0$,  $\mathrm{Res}(z)$ cannot be a polynomial in $t'$. If $\theta_1\rho_1\neq 0$, to have $\mathrm{Res}(z)$ polynomial in $t'$, we need that $\frac{-\theta_2}{\theta_1}$, respectively $\frac{-\rho_2}{\rho_1}$, is a root of order $d_0-d_1$, respectively $d_2-d_0$, of the numerator. So we would have $\frac{-\theta_2}{\theta_1}$ root of $\rho_1t'+\rho_2$, and $\frac{-\rho_2}{\rho_1}$ root of $\theta_1t'+\theta_2$, contradicting the fact that $\rho_1\theta_2-\rho_2\theta_1\neq 0$.

\noindent(2)(b) With the preceding notations, when $\min(B_f)=\max(B_f)=b_1$ (i.e. when the Newton polygon of $f$ has only one vertex $s_1=s_2$), necessarily $b_0\neq b_1$. Indeed, if not, we would, we would have  $\mathrm{Res}(z)=\mathrm{Res}(\lambda_{a_1,b_1})\in K$, which would contradict the fact that  $\mathrm{Res}(z)$ is transcendental over $K$, and consequently that $R_v$ is a dicritical divisor of $f$.
%On remarque alors, à l'instar de ??Abh-luengo, que $t$ est un autre générateur de $K_v=K(\tau)$ sur $K$ si et seulement si il existe $\phi_1,\phi_2,\psi_1,\psi_2\in K$ tels que $\tau:=\frac{\phi_1t+\phi_2}{\psi_1t+\psi_2}$. A fortiori, $\mathrm{Res}(z)$ ne peut être polynomial en $t$.
%
%
%(a) On a $s_0\notin ]s_1,s_2[$ si et seulement si les points correspondants $\tilde{s}_i=As_i$ vérifient $\tilde{s}_0\notin]\tilde{s}_1,\tilde{s}_2[$. Cela équivaut à, soit $l_0\leq l_1=\min\{l=a_2a+b_2b\ ;\ a...\}$ Supposons par exemple que $\tilde{s}_0\in]\infty,\tilde{s}_1[\subset A(D)$, ce qui implique que . On écrit le changement de coordonnées du dernier éclatement $(x_\nu,y_\nu)\mapsto (x_\nu, \frac{y_\nu}{x_\nu})$ avec $\frac{y_\nu}{x_\nu}\in R\setminus\mathfrak{m}$. En posant $t:=\mathrm{Res}\left(\frac{y}{x}\right)$, on calcule :
%$$\begin{array}{lcl}
%\mathrm{Res}(z)&=&\displaystyle\frac{\sum_{c+d= c_0+d_0}
%\mathrm{Res}(\lambda_{a,b})t^d}{t^{d_0}}\\
%&=&\sum_{c+d= c_0+d_0}
%\mathrm{Res}(\lambda_{a,b})t^{d-d_0}
%\end{array}, \ \textrm{ avec } (c,d)=(a_1 a+b_1b,a_2 a+b_2b).$$ 
%Puisque $s_0\in (\infty,s_1]\subset D$,
%Le dernier changement de coordonnées peut s'écrire $(x_\nu,y_\nu)\mapsto (x_\nu, \frac{y_\nu}{x_\nu})$ avec $\frac{y_\nu}{x_\nu}\in R\setminus\mathfrak{m}$, ou $(x_\nu,y_\nu)\mapsto ( \frac{x_\nu}{y_\nu},y_\nu)$ avec $\frac{x_\nu}{y_\nu}\in R\setminus\mathfrak{m}$. Ainsi, si on pose $\tau:=\mathrm{Res}(\frac{y}{x})$, on obtient dans le premier cas :
%$$\mathrm{Res}(z)=\sum_{c+d= c_0+d_0}
%\mathrm{Res}(\lambda_{a,b})t^b,\ \ \ t:=\mathrm{Res}\left(\frac{y}{x}\right).$$ 
\end{pf}

\section{The polynomial case.}\label{sect:jacob}
In this section, we resume the notion of dicritical divisor introduced in \cite[Section (6.2)]{abh-luengo_dicrit-div} in the case of the ring  $k[x,y]$ of bivariate polynomials over a field $k$. This notion is adapted to the Jacobian problem in dimension 2.
\begin{definition}\label{defi:dicrit-poly}
Let $f\in k[x,y]\setminus k$. We call  \textbf{ dicritical divisor} of $f$ any discrete valuation ring  $R_v$ of $k(x,y)$ such that  $k[x,y]\nsubseteq R_v$ and $k(f)\subset R_v$ with $\mathrm{Res}(f)$ transcendental over $k$.
\end{definition}
A polynomial map $f\in k[x,y]$ is defined everywhere but at infinity, where it becomes a rational function. Let  $F(X:Y:Z)=Z^mf(\frac{X}{Z},\frac{X}{Z})$, $m=\deg(f)$, be the homogenized of $f$ on $\mathds{P}_2(k)$. This function has  points of indetermination $\{F(X:Y:Z)=Z=0\}$, which are the points at infinity of the curve defined by $f$. The center of $v$ is in $\mathrm{Spec} (k[\phi,\psi])$ with $(\phi,\psi)=(\frac{1}{x},\frac{y}{x})$ if $x\notin R_v$ (open set $X\neq 0$ of $\mathds{P}_2(k)$) or $(\phi,\psi)=(\frac{1}{y},\frac{x}{y})$ if $x\in R_v$ (open set $Y\neq 0$ of  $\mathds{P}_2(k)$). For instance, in the first case, we obtain:
\begin{center}
$\begin{array}{lcl}
f(x,y)=\tilde{f}(\phi,\psi)&=&\sum_{a,b}\lambda_{a,b}x^ay^b\\
&=&x^m\sum_{a,b}\lambda_{a,b}(\frac{1}{x})^{m-(a+b)}(\frac{y}{x})^b\\
&=&\displaystyle\frac{\sum_{a,b}\lambda_{a,b}\phi^{m-(a+b)}\psi^b}{\phi^m}.
\end{array}$
\end{center}
\noindent Thus $f$ defines a rational function $\tilde{f}(\phi,\psi)\in QF(R)=k(x,y)$ for which $R_v$ is a dicritical divisor and such that $\phi^m\tilde{f}\in R$. In other words, a dicritical divisor of $f$ in the sense of (\ref{defi:dicrit-poly}) corresponds to a dicritical divisor of  $\tilde{f}$ in the sense of (\ref{defi:dicrit}) where $R$ is the local ring at a point at infinity of $f(x,y)=0$. Moreover, at these points the hypotheses of  (\ref{mainresult}) hold for $\tilde{f}$. As in the preceding section, we denote $K_v=k'(t)$ with $k'$ relative algebraic closure of $k$ in $K_v$ and $t$ transcendental over $k$. We deduce that:
\begin{Corollary}
Let $R_v$ a dicritical divisor of $f\in k[x,y]\setminus k$ in the sense of (\ref{defi:dicrit-poly}). Then the element $t$ of (\ref{rem:eclts}) can be chosen so that  $\mathrm{Res}(f)\in k'[t]$.
\end{Corollary}
This result can be seen as a complement to the study of the dicritical divisors at infinity for polynomials in two complex   variables \cite{fourrier-l_topo-poly-2var-infini}. The existence of these divisors in the general case (\ref{defi:dicrit-poly}) is not obvious. We leave to the reader the pleasure of reading the masterful argument (I) of \cite[Section (6.2)]{abh-luengo_dicrit-div}.
\def\cprime{$'$}
\providecommand{\bysame}{\leavevmode\hbox to3em{\hrulefill}\thinspace}
\providecommand{\MR}{\relax\ifhmode\unskip\space\fi MR }
\providecommand{\MRhref}[2]{%
  \href{http://www.ams.org/mathscinet-getitem?mr=#1}{#2}
}
\providecommand{\href}[2]{#2}


\begin{thebibliography}{Abh56}

\bibitem[Abh56]{abh:val-cent-local-domain}
S.~S. Abhyankar, \emph{On the valuations centered in a local domain}, Amer. J.
  Math. \textbf{78} (1956), 321--348.

\bibitem[AL11]{abh-luengo_dicrit-div}
S.~S. Abhyankar and I.~Luengo, \emph{Algebraic theory of dicritical divisors},
  to appear in Amer. Jour. Math (2011).

\bibitem[Dul06]{dulac_pts-dicrit}
H.~Dulac, \emph{Sur les points dicritiques}, J. Math. Pures Appl. (6)
  \textbf{71} (1906), no.~2, 381--402.

\bibitem[Fou96]{fourrier-l_topo-poly-2var-infini}
L.~Fourrier, \emph{Topologie d'un polyn\^ome de deux variables complexes au
  voisinage de l'infini}, Ann. Inst. Fourier (Grenoble) \textbf{46} (1996),
  no.~3, 645--687. \MR{1411724 (98e:32066)}

\bibitem[Har77]{hartshorne_alg-geom}
R.~Hartshorne, \emph{Algebraic geometry}, Springer-Verlag, New York, 1977,
  Graduate Texts in Mathematics, No. 52.

\bibitem[Hir67]{hironaka:char}
H.~Hironaka, \emph{Characteristic polyhedra of singularities}, J. Math. Kyoto
  Univ. \textbf{7} (1967), 251--293.

\bibitem[MM80]{mattei-moussu:holonomie-int-1ere}
J.-F. Mattei and R.~Moussu, \emph{Holonomie et int\'egrales premi\`eres}, Ann.
  Sci. \'Ecole Norm. Sup. (4) \textbf{13} (1980), no.~4, 469--523.

\bibitem[Sei68]{seidenberg:reduction}
A.~Seidenberg, \emph{Reduction of singularities of the differential equation
  {$A\,dy=B\,dx$}}, Amer. J. Math. \textbf{90} (1968), 248--269.

\bibitem[Spi90]{spivakovsky_funct-fields-surf}
M.~Spivakovsky, \emph{Valuations in function fields of surfaces}, Amer. J.
  Math. \textbf{112} (1990), no.~1, 107--156.

\bibitem[Tei03]{teissier:valuations}
B.~Teissier, \emph{Valuations, deformations, and toric geometry}, Valuation
  theory and its applications, {V}ol. {II} ({S}askatoon, {SK}, 1999), Fields
  Inst. Commun., vol.~33, Amer. Math. Soc., Providence, RI, 2003, pp.~361--459.

\bibitem[TW94]{le-weber_jacobian-dim2}
L\^e~Dung Tr\'ang and C.~Weber, \emph{{A geometrical approach to the Jacobian
  conjecture for $n=2$}}, Kodai Math. J. \textbf{17} (1994), no.~3, 374--381.

\end{thebibliography}
\end{document}